\documentclass[11pt]{article}
\usepackage[T1]{fontenc}
\usepackage{lmodern}
\usepackage{amsmath,amssymb,amsthm,mathtools}
\usepackage{microtype}
\usepackage{xurl}
\usepackage{xcolor}
\usepackage{enumitem}
\usepackage[margin=1in]{geometry}
\usepackage{hyperref}
\hypersetup{
  colorlinks=true,
  linkcolor=blue!55!black,
  citecolor=blue!55!black,
  urlcolor=blue!55!black,
  pdftitle={From Internal to External: Classical Models of ZF+PP+not-AC -- Retraction of the Proposed Construction},
  pdfauthor={Frank Gilson}
}

\newcommand{\ZF}{\mathsf{ZF}}
\newcommand{\AC}{\mathsf{AC}}
\newcommand{\PP}{\mathsf{PP}}

\newcommand{\Sh}{\mathsf{Sh}}
\newcommand{\B}{\mathbb B}
\newcommand{\Pp}{\mathbb P}
\newcommand{\NN}{\mathbb N}
\newcommand{\onto}{\twoheadrightarrow}
\newcommand{\into}{\hookrightarrow}
\newcommand{\id}{\operatorname{id}}
\DeclareMathOperator{\Fn}{Fn}
\DeclareMathOperator{\Fix}{Fix}
\DeclareMathOperator{\supp}{supp}

\theoremstyle{plain}
\newtheorem{proposition}{Proposition}[section]
\newtheorem{lemma}[proposition]{Lemma}

\title{From Internal to External: Classical Models of
\texorpdfstring{$\ZF+\PP+\neg\AC$}{ZF+PP+not-AC}\\
\large Retraction of the Proposed Construction}
\author{Frank Gilson}
\date{September 10, 2026}

\begin{document}
\maketitle

\begin{abstract}
An earlier form of this manuscript claimed to obtain a classical symmetric
model of \(\ZF+\PP+\neg\AC\) by two routes: a Boolean-valued
externalization of an internal sheaf-topos construction, and a direct
symmetric extension based on \(\Fn(\NN\times H,2)\). Those claims are
withdrawn.

The decisive contradiction is already present in the stated
Local-to-Global Embedding Principle. For every hereditarily symmetric
surjection \(f:A\onto B\), that theorem purports to construct a function
\(s:B\to A\) satisfying \(f\circ s=\id_B\). Thus it asserts that every
surjection has a right inverse, which is equivalent over \(\ZF\) to the
Axiom of Choice and is stronger than the injection required by the
Partition Principle. It therefore cannot prove \(\PP+\neg\AC\).

The supporting localization and gluing arguments fail independently: a
condition cannot be strengthened while deleting coordinates outside a
fixed support; the fixed support is replaced in the proof by
element-dependent supports; the proposed coherent antichains cannot exist
below arbitrary Boolean values; minimality in an arbitrary ground-model
well-order does not imply equivariance; right inverses need not agree on
overlaps; and countably many finite supports need not combine into a
finite support. The route-unification theorem and the proposed witness to
\(\neg\AC\) have further defects. This replacement retains only standard
background facts about symmetric extensions and records the obstructions.
It asserts no model of \(\ZF+\PP+\neg\AC\) and no consistency or
independence result.
\end{abstract}



\section*{Correction and current status}

The earlier manuscript advertised:
\begin{enumerate}[label=(\roman*)]
  \item two equivalent constructions of a symmetric model \(N\);
  \item a Local-to-Global Embedding Principle for hereditarily symmetric
        names;
  \item \(N\models\PP\);
  \item \(N\models\neg\AC\); and therefore
  \item \(N\models\ZF+\PP+\neg\AC\), with corresponding
        relative-consistency and nonimplication claims.
\end{enumerate}

\noindent\textbf{All five conclusions are withdrawn.}
The Local-to-Global Embedding Principle is incompatible with the claimed
failure of choice, and neither route supplies a valid independent proof.
The former main theorem, the route-unification theorem, all corollaries
depending on them, and the purported consistency conclusion must not be
cited as established results.

The standard definitions of forcing, symmetric names, normal filters of
subgroups, and Boolean-valued semantics remain useful background. The
standard symmetric-extension theorem also remains available: a genuine
symmetric system yields a transitive model of \(\ZF\). None of those
background facts establishes \(\PP\), failure of \(\AC\), or the
advertised separation.

\section{The central strength mismatch}

Recall that the Partition Principle is
\[
  \PP:\qquad
  \forall A\,\forall B\,
  \bigl((\exists f:A\onto B)\longrightarrow
        (\exists j:B\into A)\bigr).
\]
The injection \(j\) need not select a point in the corresponding fiber of
\(f\). Requiring it to do so changes the principle completely.

\begin{proposition}\label{prop:sections-ac}
Over \(\ZF\), the following are equivalent:
\begin{enumerate}[label=(\arabic*)]
  \item \(\AC\);
  \item every surjection has a right inverse.
\end{enumerate}
\end{proposition}

\begin{proof}
Assume \(\AC\), and let \(f:A\onto B\). Choosing one element from every
nonempty fiber \(f^{-1}(\{b\})\) gives \(s:B\to A\) with
\(f\circ s=\id_B\).

Conversely, let \(\langle X_i:i\in I\rangle\) be a set-indexed family of
nonempty sets and put
\[
  E=\{\langle i,x\rangle:i\in I\text{ and }x\in X_i\}.
\]
The projection \(E\onto I\) is surjective. A right inverse has the form
\(i\mapsto\langle i,c(i)\rangle\), where \(c(i)\in X_i\), and hence
supplies a choice function.
\end{proof}

The former theorem labelled \emph{LEP for symmetric names} asserted that,
for every hereditarily symmetric surjection name
\(\dot f:\dot A\onto\dot B\), there is a hereditarily symmetric
\(\dot s:\dot B\to\dot A\) with
\[
  \Vdash \dot f\circ\dot s=\id_{\dot B}.
\]
The later theorem labelled \emph{PP in \(N\)} repeated this right-inverse
conclusion. By Proposition~\ref{prop:sections-ac}, the asserted theorem
would give \(N\models\AC\), not merely \(N\models\PP\). It therefore
contradicts the advertised \(N\models\neg\AC\).

The same problem reappears in the former uniform symmetric collection
lemma. From
\[
  \Vdash \forall x\in\dot X\ \exists y\ \varphi(x,y)
\]
it purports to construct one symmetric function choosing a witness for
every \(x\). Applied to membership in a family of nonempty sets, that is a
choice principle rather than an instance of Collection. Pointwise
existence of witnesses is not an internal family of witnesses.

\section{Why fixed-support localization fails}

Let
\[
  \Pp=\Fn(\NN\times H,2)
\]
be ordered so that stronger conditions extend weaker ones, and let
\(\supp(p)\) be the finite set of rows used by \(p\).

The former uniform-localization condition required, for a fixed finite
support \(S\), that arbitrary relevant conditions \(p\) have extensions
\(q\leq p\) satisfying \(\supp(q)\subseteq S\). This is impossible as
soon as \(p\) mentions a row outside \(S\), because
\[
  q\leq p\quad\Longrightarrow\quad
  \operatorname{dom}(q)\supseteq\operatorname{dom}(p)
  \quad\Longrightarrow\quad
  \supp(q)\supseteq\supp(p).
\]
No thinning argument can remove coordinates from a stronger condition.
In particular, conditions supported by \(S\) are not dense in \(\Pp\).

There are three further support failures.
\begin{enumerate}[label=(\arabic*)]
\item The fixed-support local-section lemma begins with one finite support
\(S\) intended to support a section on all of \(B\), but its proof fixes
\(b\in B\) and replaces \(S\) by a support \(S_b\) of that element. It
therefore produces, at best, element-dependent local data and no single
function with one finite support.

\item A support of \(\dot b\) is treated as if it also supported
\(\dot f\). Hereditary symmetry of \(\dot f\) supplies some finite
support, but not necessarily the chosen support of \(\dot b\). Enlarging
the support repairs that bookkeeping error only for one pair
\((\dot f,\dot b)\); it does not make the family uniform.

\item From the existence of an extension \(q\leq p\) forcing an
existential statement, the proof concludes that \(p\) forces the
statement. One stronger condition witnessing compatibility does not show
that every generic containing \(p\) satisfies the assertion.
\end{enumerate}

The forcing automorphism lemma has the form
\[
  p\Vdash\varphi(\dot x)
  \quad\Longleftrightarrow\quad
  \pi p\Vdash\varphi(\pi\dot x).
\]
It transports the condition as well as the name. To compare both claims
below the same \(p\), one still needs a hypothesis such as \(\pi p=p\),
or a compatibility argument producing a common extension. Support of a
name alone does not give same-generic or same-condition equivariance.
This defect occurs in the localization arguments of both routes.

\section{Failure of the antichain and least-witness arguments}

Route~I introduced an \(S\)-coherent antichain below a Boolean value
\(b\), requiring every antichain element \(u\) to be fixed by
\(\Fix(S)\), and claimed that such an antichain could be maximal below
an arbitrary \(b\).

\begin{lemma}\label{lem:antichain-obstruction}
Let a group act by complete automorphisms on a complete Boolean algebra
\(\B\). If every member of \(A\subseteq\B\) is fixed by \(\Fix(S)\), then
\(\bigvee A\) is fixed by \(\Fix(S)\).
\end{lemma}

\begin{proof}
For every \(\pi\in\Fix(S)\),
\[
  \pi\!\left(\bigvee A\right)
  =\bigvee_{u\in A}\pi(u)
  =\bigvee A.
\]
\end{proof}

Thus an \(S\)-coherent maximal antichain below \(b\) whose join is \(b\)
can exist only when \(b\) itself is \(S\)-invariant. The former lemma
asserted it for arbitrary Boolean values. Extending an \(S\)-cylinder
does not make a Boolean value fixed by every automorphism in
\(\Fix(S)\).

The selector construction then chose the least witness under a fixed
ground-model well-order and declared the choice equivariant. An arbitrary
well-order of the forcing objects need not be preserved by the
automorphism group:
\[
  x\mathrel{\triangleleft}y
  \quad\not\Longrightarrow\quad
  \pi x\mathrel{\triangleleft}\pi y.
\]
Applying \(\pi\) to a least witness therefore need not yield the least
witness for the transported problem. Ground-model definability of the
well-order does not supply invariance under permutations of its elements.
The mixing, selector, and Boolean-valued localization constructions using
this step are unproved.

\section{Failure of coherence and gluing}

Even if local right inverses were available, their asserted coherence
does not follow.

\begin{proposition}\label{prop:right-inverses}
Two right inverses of the same surjection need not agree.
\end{proposition}

\begin{proof}
Let \(f:\{0,1\}\onto\{\ast\}\) be the unique map. The maps
\(s_0(\ast)=0\) and \(s_1(\ast)=1\) both satisfy
\(f\circ s_i=\id_{\{\ast\}}\), but \(s_0\neq s_1\).
\end{proof}

Separativity of the forcing does not alter this elementary fact. Hence
local sections do not agree on overlaps merely because their
postcomposition with \(f\) is the identity.

The support-directed gluing lemma has two additional defects. If
\(S\subseteq T\), then
\[
  \Fix(T)\subseteq\Fix(S),
\]
not the reverse inclusion used in the proof. More importantly, gluing
through a countable increasing family of finite supports may produce a
name supported by their countable union. A countable union of finite
supports need not be finite. The normal filter generated by finite
stabilizers does not thereby contain the stabilizer of that countable
union.

Accordingly, countably many locally symmetric choices do not assemble
into one hereditarily symmetric function by the proposed argument.

\section{Failure of Route I}

Route~I depended on the earlier manuscript
\emph{Internal models of \(\ZF+\PP+\neg\AC\) via LEP in Boolean topoi}
\cite{InternalPPTopos2025}. The current replacement of that manuscript
withdraws its proposed internal model, local embedding theorem, witness
to failure of choice, and Booleanization upgrade. Thus the input that
Route~I was supposed to externalize is no longer available.

There is also an independent error at the start of Route~I. For a free
action of a finite group \(H\) on a space \(X\), the action groupoid is
Morita equivalent to the quotient presentation, so the associated topos
is equivalent to \(\Sh(X/H)\). A sheaf topos on a Cantor-type quotient is
not generally Boolean: arbitrary open sets need not have open
complements. The Boolean algebra of clopen subsets is generally not
complete. Passing instead to the complete Boolean algebra of regular
open sets constructs a Booleanization; it does not identify the original
sheaf topos with that Boolean-valued presentation.

More generally, a Boolean Grothendieck topos need not be the topos of
sheaves on a complete Boolean algebra unless the required localic
hypotheses are supplied. The former Boolean-presentation proposition was
therefore too broad. These issues prevent Route~I from serving as an
independent transfer argument.

\section{Failure of the route-unification theorem}

The proposed translation of Boolean names to poset names used the
equivalence
\[
 b\Vdash\varphi(\dot y)
 \quad\Longleftrightarrow\quad
 \exists p\in A(b)\ 
 p\Vdash\varphi(T(\dot y)),
\]
where \(A(b)\) is a maximal antichain below \(b\). The existential
right-hand side is too weak. For \(b\) to force the formula, all relevant
pieces below \(b\), equivalently a dense family below \(b\), must force
it; one favorable antichain member is insufficient.

The claimed equivariance of the translation uses the non-invariant
least-witness selection discussed above. The functoriality argument also
asserts that if \(i:\Pp\to\B^+\) is dense and
\(\rho:\B\to\mathbb C\) is a complete homomorphism, then
\(\rho\circ i\) is dense in \(\mathbb C\). It is only necessarily dense
in \(\rho[\B]\), which may be a proper complete subalgebra of
\(\mathbb C\). A complete Boolean homomorphism also does not by itself
transport the chosen automorphism group, normal filter, and support
calculus.

The alleged canonical equivalence of the two symmetric models is
therefore not established.

\section{Failure of the proposed witness to
\texorpdfstring{\(\neg\AC\)}{not-AC}}

The direct-forcing route defined a coordinate real \(\dot r_c\) from one
forcing bit at a coordinate \(c\), using two ground-model numerical labels
\(\alpha(c)\) and \(\beta(c)\). It then claimed that a row transposition
taking \(c_0\) to \(c_1\) swaps \(\dot r_{c_0}\) and
\(\dot r_{c_1}\).

Forcing automorphisms move conditions and coordinate occurrences in
names, but fix canonical names for ground-model natural numbers. Thus the
automorphism moves the bit controlling the definition while leaving the
labels \(\alpha(c_0)\) and \(\beta(c_0)\) unchanged. It does not thereby
turn the name into the differently labelled \(\dot r_{c_1}\). Each
evaluated \(\dot r_c\) in the displayed construction is merely one of
two ground-model singleton reals determined by those labels. The asserted
transposition calculation is false, and the resulting family does not
establish non-well-orderability.

The forcing \(\Fn(\NN\times H,2)\) also provides only finitely many bits
in each row because \(H\) is finite. The manuscript did not supply a
different symmetric object with the independent coordinate structure
needed for a standard Cohen-style failure-of-choice argument. Hence
\(N\models\neg\AC\) is not proved.

\section{Status of the former main theorem}

The failures above have different logical status, but converge on the
same public conclusion.

\begin{center}
\begin{tabular}{p{0.30\textwidth}p{0.59\textwidth}}
\hline
\textbf{Former component} & \textbf{Current status}\\
\hline
LEP right-inverse theorem &
False as a route to \(\PP+\neg\AC\): its conclusion is equivalent to
\(\AC\).\\[2pt]
Uniform localization &
False for the stated finite-condition forcing.\\[2pt]
Coherent-antichain lemma &
False for arbitrary Boolean values.\\[2pt]
Least-witness equivariance &
Not established; an arbitrary well-order is not automorphism-invariant.\\[2pt]
Local-section coherence &
False; right inverses need not agree.\\[2pt]
Support-directed gluing &
Proof failure; subgroup inclusion is reversed and finite support is lost.\\[2pt]
Route I &
Dependency failure and independent presentation errors.\\[2pt]
Route unification &
Proof failure.\\[2pt]
Coordinate-real argument &
False automorphism calculation; \(\neg\AC\) is unproved.\\
\hline
\end{tabular}
\end{center}

No weakening of the advertised theorem is established by the remaining
argument. Replacing right inverse by injection would remove the immediate
implication of \(\AC\), but every construction of the injection would
then have to be rebuilt; the existing fiber-selection proof cannot
simply be relabelled.

\section{Conclusion}

The proposed internal-to-external construction does not establish a
classical model of
\[
  \ZF+\PP+\neg\AC.
\]
Its central theorem aims at the wrong target and would imply \(\AC\);
the localization and gluing machinery does not construct the asserted
symmetric selectors; Route~I has no valid internal model to externalize;
the route-unification theorem fails; and the direct argument does not
prove \(\neg\AC\).

Accordingly, the original construction and every result depending on it
are superseded by this correction. This paper makes no relative
consistency claim and no claim that \(\PP\) fails to imply \(\AC\) over
\(\ZF\). The question is not resolved by the arguments presented here.

\end{document}